\documentstyle[12pt,psfig,amssymb,amsmath]{article}

\newcommand{\be}{\begin{equation}}
\newcommand{\ee}{\end{equation}}
\newcommand{\bea}{\begin{eqnarray}}
\newcommand{\eea}{\end{eqnarray}}

\def\1#1{^{(#1)}}
\def\la{\langle}
\def\ra{\rangle}

\begin{document}
\title{Power Sums Related to Semigroups ${\sl S}\left(d_1,d_2,d_3\right)$}
\author{Leonid G. Fel${}^{\dag}$ and Boris Y. Rubinstein${}^{\ddag}$\\
\\
${}^{\dag}$Department of Civil and Environmental Engineering,\\
Technion, Haifa 32000, Israel\\
and \\
${}^{\ddag}$Department of Mathematics, University of California, Davis,
\\One Shields Dr., Davis, CA 95616, U.S.A.}
\date{\today}

\maketitle
\def\be{\begin{equation}}
\def\ee{\end{equation}}
\def\bea{\begin{eqnarray}}
\def\eea{\end{eqnarray}}
\def\p{\prime}
%\vspace{-1cm}
\begin{abstract}
The explicit formulas for the sums of positive powers of the
integers $s_i$ unrepresentable by the triple of integers
$d_1,d_2,d_3\in {\mathbb N},\;\gcd(d_1,d_2,d_3)=1$, are derived.
\vskip .1 cm 
\noindent 
\begin{tabbing}
{\sf Key words}:\hspace{.1in}  \=Non--symmetric and symmetric semigroups, 
Power sums.
\end{tabbing}    
\vskip .1 cm
{\sf 2000 Math. Subject Classification}: Primary - 11P81; Secondary - 20M99
\end{abstract}
\vspace{.1cm}

\newpage
Let ${\sf S}\left(d_1,\ldots,d_m\right)$ be the semigroup generated by a set
of integers $\{d_1,\ldots,d_m\}$ such that
\begin{equation}
1<d_1<\ldots <d_m\;,\;\;\;\gcd(d_1,\ldots,d_m)=1\;.
\label{definn0}
\end{equation}
For short we denote the tuple $(d_1,\ldots,d_m)$ by ${\bf d}^m$ and consider 
the generating function $\Phi\left({\bf d}^m;z\right)$ 
\begin{eqnarray}
\Phi\left({\bf d}^m;z\right)=\sum_{s_i\;\in\;\Delta\left({\bf d}^m\right)} 
z^{s_i}=z+z^{s_2}+\ldots+z^{s_{G\left({\bf d}^m\right)}}\;,\label{def1}
\end{eqnarray}
for the set  $\Delta\left({\bf d}^m\right)$ of the integers $s$ 
which are unrepresentable as $s=\sum_{i=1}^mx_id_i,\;x_i\in {\mathbb 
N}\cup\{0\}$
\begin{eqnarray}
\Delta\left({\bf d}^m\right)=\{s_1,s_2,\ldots,s_{G\left({\bf d}^m\right)}\}\;,
\;\;\;s_1=1\;.
\label{def2}
\end{eqnarray}
The integer $G\left({\bf d}^m\right)$ is known as the genus for semigroup 
${\sf S}\left({\bf d}^m\right)$.  Recall the relation of $\Phi\left({\bf d}^m
;z\right)$ with the Hilbert series $H({\bf d}^m;z)$ of a graded monomial 
subring ${\sf k}\left[z^{d_1},\ldots,z^{d_m}\right]$ \cite{stan96}
\begin{equation}
H({\bf d}^m;z)+\Phi({\bf d}^m;z)=\frac{1}{1-z}\;,\;\;\;\;\mbox{where}\;\;\;
H({\bf d}^m;z)=\sum_{s\;\in\;{\sf S}\left({\bf d}^m\right)} z^s=
\frac{Q({\bf d}^m;z)}{\prod_{j=1}^m \left(1-z^{d_j}\right)}\;,
\label{gen0}
\end{equation}
and $Q({\bf d}^m;z)$ is a polynomial in $z$. The calculation of the power 
sums  
\begin{equation}
g_n\left({\bf d}^m\right)=\sum_{s_i\;\in\;\Delta\left({\bf d}^m\right)}
s_i^n\label{gen1}
\end{equation}
was performed in \cite{rods78} for $m=2$
\begin{eqnarray}
g_n\left({\bf d}^2\right)=\frac{1}{(n+1)(n+2)}\sum_{k=0}^{n+1}\sum_{l=0}^{n+1-k}
{n+2\choose k}{n+2-k\choose l}B_kB_ld_1^{n+1-k}d_2^{n+1-l}-\frac{B_{n+1}}{n+1}
\;,\label{gen1a}
\end{eqnarray}
where $B_k$ are the Bernoulli numbers. The formula (\ref{gen1a}) generalizes 
the known Sylvester's expression \cite{sylv84} for $G\left({\bf d}^2\right)=
g_0\left({\bf d}^2\right)$ and further results for $n=1$ \cite{jau93} and 
$n=2,3$ \cite{fel04}
\begin{eqnarray}
g_0\left({\bf d}^2\right) & = & \frac{1}{2}(d_1-1)(d_2-1)\;,
\nonumber \\
g_1\left({\bf d}^2\right) & = &
\frac{g_0\left({\bf d}^2\right)}{6}\;
(2d_1d_2-d_1-d_2-1)\;,
\nonumber \\
g_2\left({\bf d}^2\right) & = &
\frac{g_0\left({\bf d}^2\right)}{6}\;
d_1d_2(d_1d_2-d_1-d_2)\;,
\nonumber \\
g_3\left({\bf d}^2\right) & = &
\frac{g_0\left({\bf d}^2\right)}{60}\;
\left[(1+d_1)(1+d_2)\left(1+d_1^2+d_2^2+6d_1^2d_2^2\right)-
15d_1d_2 (d_1+d_2)\right]\;.\nonumber
\end{eqnarray}
As for higher dimensions, $m\geq 3$, the first two sums, $g_0\left({\bf d}^3
\right)$ and $g_1\left({\bf d}^3\right)$, were found in \cite{fel04}. The use 
was made of the explicit expression for the Hilbert series $H({\bf d}^3;z)$ of
a graded subring for semigroups ${\sf S}\left({\bf d}^3\right)$ which was
recently established \cite{fel04}.

In this paper we derive the formula for the power sum $g_n\left({\bf d}^3
\right)$ related to the symmetric and non--symmetric semigroups ${\sl S}\left(
{\bf d}^3\right)$. This will be done by applying to the relation (\ref{gen0}) 
an approach \cite{rods78} exploiting the properties of the generating 
function of the Bernoulli numbers.
%\begin{eqnarray}
%\frac{x}{e^x-1}=\sum_{k=0}^{\infty}B_k\frac{x^k}{k!}\;.
%\label{bern1}
%\end{eqnarray}

Following Johnson \cite{john60} recall the definition of {\em the minimal 
relations} for given $d_1,d_2,d_3$ 
\begin{eqnarray}
&&a_{11}d_1=a_{12}d_2+a_{13}d_3\;,\;\;\;a_{22}d_2=a_{21}d_1+a_{23}d_3\;,\;\;\;
a_{33}d_3=a_{31}d_1+a_{32}d_2\;,\;\;\;\;\mbox{where}\;\;\;\;\;\label{herznon1a}\\
&&a_{11}=\min\left\{v_{11}\;\bracevert\;v_{11}\geq 2,\;v_{11}d_1=v_{12}d_2+
v_{13}d_3,\;v_{12},v_{13}\in {\mathbb N}\cup\{0\}\right\}\;,\nonumber\\
&&a_{22}=\min\left\{v_{22}\;\bracevert\;v_{22}\geq 2,\;v_{22}d_2=v_{21}d_1+
v_{23}d_3,\;v_{21},v_{23}\in {\mathbb N}\cup\{0\}\right\}\;,\label{joh1}\\
&&a_{33}=\min\left\{v_{33}\;\bracevert\;v_{33}\geq 2,\;v_{33}d_3=v_{31}d_1+ 
v_{32}d_2,\;v_{31},v_{32}\in {\mathbb N}\cup\{0\}\right\}\;.\nonumber
\end{eqnarray}
The auxiliary invariants $a_{ij}$ are uniquely defined by (\ref{joh1}) and 
\begin{eqnarray}
\gcd(a_{11},a_{12},a_{13})=1\;,\;\;\;\gcd(a_{21},a_{22},a_{23})=1\;,\;\;\;
\gcd(a_{31},a_{32},a_{33})=1\;.\label{joh2}
\end{eqnarray}
According to \cite{fel04} the numerator of the Hilbert series for 
non--symmetric semigroups ${\sf S}\left({\bf d}^3\right)$ reads
\begin{eqnarray}
Q({\bf d}^3;z)&=&
1-\sum_{i=1}^3 z^{a_{ii}d_i}+
z^{1/2\left[\langle{\bf a},
{\bf d}\rangle-J\left({\bf d}^3\right)\right]}
+z^{1/2\left[\langle{\bf a},{\bf
d}\rangle+J\left({\bf d}^3\right)\right]}\;,\;\;\;\mbox{where}\label{nomn1}\\
J^2\left({\bf d}^3\right)&=&\langle{\bf a},{\bf d}\rangle^2-4\sum_{i>j}^3
a_{ii}a_{jj}d_id_j+4d_1d_2d_3\;,\;\;\;\;\langle{\bf a},{\bf d}\rangle=
\sum_{i=1}^3a_{ii}d_i\;.\label{nomn2}
\end{eqnarray}
The case of symmetric semigroups ${\sf S}\left({\bf d}^3\right)$ is much
 simpler. 
Here two off--diagonal elements in one column of the matrix $a_{ij}$ vanish, 
e.g. $a_{13}=a_{23}=0$ that results in $a_{11}d_1=a_{22}d_2$. The numerator of 
the Hilbert series for symmetric semigroups ${\sf S}\left({\bf d}^3\right)$ 
with above symmetry is given by
\begin{eqnarray}
Q_s({\bf d}^3;z)=(1-z^{a_{22}d_2})(1-z^{a_{33}d_3})\;.\label{nomn3}
\end{eqnarray}  
Move on to calculation of the power sums $g_n\left({\bf d}^3\right)$. 
First, denote $z=e^t$ and represent (\ref{gen0}) as follows
\begin{eqnarray}
\sum_{s_i\;\in\;\Delta\left({\bf d}^3\right)}e^{s_it}=\frac{1}{1-e^t}-
\frac{Q({\bf d}^3;e^t)}{\left(1-e^{d_1t}\right)\left(1-e^{d_2t}\right)
\left(1-e^{d_3t}\right)}\;,\label{gen3}
\end{eqnarray}
and apply the sequence of identities
\begin{eqnarray}
\sum_{s_i\;\in\;\Delta\left({\bf d}^3\right)}e^{s_it}=
\sum_{s_i\;\in\;\Delta\left({\bf d}^3\right)}
\sum_{k=0}^{\infty}\frac{s_i^nt^n}{n!}=\sum_{k=0}^{\infty}\frac{t^n}{n!}
\sum_{s_i\;\in\;\Delta\left({\bf d}^3\right)}s_i^n=
\sum_{n=0}^{\infty}g_n\left({\bf d}^3\right)\frac{t^n}{n!}\;.\label{gen4}
\end{eqnarray}
Thus, $g_n\left({\bf d}^3\right)$ could be found by expanding the right
hand side of (\ref{gen3}) in the power series.

In order to present the computation in the compact form we need the definition
of the Bernoulli polynomials of higher order \cite{bat53}
%\begin{eqnarray}
%\frac{t^m e^{tx} \prod_{i=1}^m d_i}{\prod_{i=1}^m (e^{d_i t}-1)} =
%\sum_{n=0}^{\infty} \frac{t^n}{n!} B_{n}^{(m)}(x|{\bf d}^m),
%\label{Bern_poly_HO_defGF}
%\end{eqnarray}
%which gives rise to the formula
\begin{eqnarray}
\frac{e^{xt}}{\prod_{i=1}^m (e^{d_i t}-1)} =
\frac{1}{\pi_d}\sum_{n=0}^{\infty}
\frac{t^{n-m}}{n!} B_{n}^{(m)}(x|{\bf d}^m)\;,\;\;\;\pi_d=\prod_{i=1}^m d_i,
\label{Bern_poly1}
\end{eqnarray}
which satisfy the recursion relation 
\begin{eqnarray}
B_{n}^{(m)}(x|{\bf d}^m) =\sum_{p=0}^n {n\choose p} d_m^p B_p
B_{n-p}^{(m-1)}(x|{\bf d}^{m-1}), \ \
B_n^{(1)}(x|d) = d^n B_n\left(\frac{x}{d}\right),
\label{Bernrec}
\end{eqnarray}
where ${\bf d}^{m-1}$ denotes the tuple $(d_1,\ldots,d_{m-1})$ and $B_n(x)$ stands for 
the regular Bernoulli polynomial, and $B_n\left(0\right)=B_n$. Because each term in 
the right hand side of (\ref{gen3}) has the form of the left hand side of 
(\ref{Bern_poly1}), it is easy to write the answer as sum of the Bernoulli 
polynomials of higher order. 

First, contribution of the term $1/(1-z)$ to $g_n\left({\bf d}^3\right)$ is 
found trivially
\begin{eqnarray}
\frac{1}{1-e^{t}} =-\sum_{n=0}^{\infty}\frac{t^{n-1}}{n!} B_{n}^{(1)}(0|1) =
-\sum_{n=0}^{\infty}\frac{t^{n-1}}{(n-1)!} \frac{B_{n}}{n},
\label{term0}
\end{eqnarray}
so that the corresponding term in $g_n\left({\bf d}^3\right)$ is 
$-B_{n+1}/(n+1)$. 
%\begin{eqnarray}
%-\frac{B_{n+1}}{n+1}.\label{term00}
%\end{eqnarray}
Consider a general term for $n=3$
\begin{eqnarray}
\frac{e^{tx}}{(1-e^{d_1 t})(1-e^{d_2 t})(1-e^{d_3 t})} =
-\frac{1}{d_1d_2d_3}\sum_{n=0}^{\infty}
\frac{t^{n-3}}{n!} B_{n}^{(3)}(x|{\bf d}^3),
\label{Bern_poly3}
\end{eqnarray}
thus its contribution to $g_n({\bf d}^3)$ reads:
\begin{eqnarray}
g(x;{\bf d}^3)=-\frac{n!}{(n+3)!\;d_1d_2d_3}B_{n+3}^{(3)}(x|{\bf d}^3).
\label{Bern_poly3a}
\end{eqnarray}
The Bernoulli polynomial of higher order $B_{n+3}^{(3)}(x|{\bf d}^3)$
 can be expanded into the triple sum over the Bernoulli numbers
\begin{eqnarray}
B_{n+3}^{(3)}(x|{\bf d}^3)=
\sum_{j=0}^{n+3}\sum_{k=0}^j\sum_{l=0}^k
{n+3\choose j}{j\choose k}{k\choose l}
d_1^{k-l}d_2^{j-k}d_3^{n+3-j} B_{k-l}B_{j-k}B_{n+3-j} x^l,
\label{B3expand}
\end{eqnarray}
which for $x=0$ reduces to the double sum
\begin{eqnarray}
B_{n+3}^{(3)}(0|{\bf d}^3)=
\sum_{j=0}^{n+3}\sum_{k=0}^j
{n+3\choose j}{j\choose k}
d_1^{k}d_2^{j-k}d_3^{n+3-j} B_{k}B_{j-k}B_{n+3-j}.
\label{B3expand0}
\end{eqnarray}
The expression for $g_n({\bf d}^3)$ for the non--symmetric semigroups has 
the form
\begin{eqnarray}
g_n^{(n)}({\bf d}^3) = 
-\frac{B_{n+1}}{n+1} + \frac{n!}{(n+3)!\;d_1d_2d_3}
\left[B_{n+3}^{(3)}(0|{\bf d}^3)-B_{n+3}^{(3)}(a_{11}d_1|{\bf d}^3)
-B_{n+3}^{(3)}(a_{22}d_2|{\bf d}^3)-
%-B_{n+3}^{(3)}(a_{33}d_3|{\bf d}^3)
\right.
\label{gnd3ns}
 \\
\left.
B_{n+3}^{(3)}(a_{33}d_3|{\bf d}^3)
+B_{n+3}^{(3)}(1/2\left[\langle{\bf a},
{\bf d}\rangle-J\left({\bf d}^3\right)\right]|{\bf d}^3)
+
B_{n+3}^{(3)}(1/2\left[\langle{\bf a},
{\bf d}\rangle+J\left({\bf d}^3\right)\right]|{\bf d}^3)
\right].
\nonumber
\end{eqnarray}
In case of the symmetric semigroups we obtain
\begin{eqnarray}
g_n^{(s)}({\bf d}^3) =
-\frac{B_{n+1}}{n+1} + \frac{n!}{(n+3)!\;d_1d_2d_3}
\left[
B_{n+3}^{(3)}(0|{\bf d}^3)
-B_{n+3}^{(3)}(a_{22}d_2|{\bf d}^3)-B_{n+3}^{(3)}(a_{33}d_3|{\bf d}^3)+
\right.
\nonumber \\
\left.
B_{n+3}^{(3)}(a_{22}d_2+a_{33}d_3|{\bf d}^3)\right].
\label{gnd3s}
\end{eqnarray}
%Below we give explicit expressions for first genera. We start from
%$g_n({\bf d}^2)$ for $n=4,5,6$
%\begin{eqnarray}
%g_4\left({\bf d}^2\right) & = &
%\frac{g_2\left({\bf d}^2\right)}{5}\;
%(2d_1d_2(d_1d_2-d_1-d_2)-d_1^2-d_2^2+d_1d_2)\;,
%\label{g4d2} \\
%g_5\left({\bf d}^2\right) & = &
%\frac{g_0\left({\bf d}^2\right)}{252}\;
%\left(
%6d_1^4d_2^4(9-5(d_1+d_2)) +
%12d_1^2d_2^2(d_1d_2(d_1^2+d_2^2)+d_1^3d_2^3+d_1^3+d_2^3)
%\right.
%\nonumber \\
%&-&
%\left.
%9d_1^2d_2^2(d_1d_2(d_1+d_2)+d_1^2+d_2^2)
%-2(1+d_1)(1+d_2)(d_1^4+d_2^4+(1+d_1^2)(1+d_2^2))
%\right),
%\label{g5d2} \\
%g_6\left({\bf d}^2\right) & = &
%\frac{g_2\left({\bf d}^2\right)}{14}\;
%\left(
%2(d_1^2+d_2^2)^2-2d_1d_2(d_1+d_2)^2+2d_1^2d_2^2+
%4 d_1d_2(d_1+d_2)^3
%\right.
%\nonumber \\
%&-&
%\left.
%d_1^2d_2^2(d_1+d_2)^2+6d_1^3d_2^3(d_1-1)(d_2-1)+
%5d_1^3d_2^3-3d_1^4d_2^4-11d_1^2d_2^2(d_1+d_2)
%\right).
%\label{g6d2}
%\end{eqnarray}
Thus, for the symmetric semigroup ${\sf S}({\bf d}^3)$ the first three sums read
\bea
2g_0^{(s)}({\bf d}^3)&=&
1-s_d+\la\widetilde{{\bf a},{\bf d}}\ra\;,\;\;\;\;\;s_d=\sum_{i=1}^3d_i\;,\;
\;\;\;\la\widetilde{{\bf a},{\bf d}}\ra=\la{\bf a},{\bf d}\ra-a_{11}d_1\;,
\label{g0d3s} \\
12g_1^{(s)}({\bf d}^3)&=&
\left(s_d-\la\widetilde{{\bf a},{\bf d}}\ra\right)\left(s_d-2\la\widetilde{
{\bf a},{\bf d}}\ra\right)+(d_1d_2+d_1d_3+d_2d_3)-\pi_d-1\;,\label{g1d3s} \\
12g_2^{(s)}({\bf d}^3)&=&(s_d-\la\widetilde{{\bf a},{\bf d}}\ra)
\left(s_d\la\widetilde{{\bf a},{\bf d}}\ra-\la\widetilde{{\bf a},{\bf d}}\ra^2
-(d_1d_2+d_1d_3+d_2d_3)+\pi_d\right)\;.\label{g2d3s} % \\
%g_3^{(s)} & = &
%\frac{1}{120}
%\left(
%1+30\la {\bf a},{\bf d}^{(4)} \ra-s_d^{(4)}+
%5\sum_{i>j}^3d_i^2d_j^2+
%5\la {\bf a},{\bf d} \ra^2(3s_d^2-s_d^{(2)})
%\right.
%\label{g3d3s} \\
%&&
%-15s_d(\la {\bf a},{\bf d} \ra^3+\la {\bf a},{\bf d}^{(3)}\ra)-
%15s_d^{(2)}\la {\bf a},{\bf d}^{(2)}\ra+
%6\la {\bf a}^{(2)},{\bf d}^{(2)}\ra^2
%-22 \pi_d^2
%\nonumber \\
%&&
%\left.
%+15\pi_d\la {\bf a},{\bf d} \ra^2-
%30\pi_d\la {\bf a},{\bf d} \ra+
%15\pi_d s_d(1-a_{22}-a_{33}+\la {\bf a},{\bf d} \ra)
%-5\pi_d(s_d^2+\sum_{i>j}^3d_id_j)
%\right)\;,
%\nonumber
%\\
%g_4^{(s)} & = &
\eea
%where
%$$
%\pi_d = \prod_{i=1}^3 d_{i}, \ \
% s_d^{(k)} = \sum_{i=1}^3 d_{i}^k, \ \
%s_d \equiv s_d^{(1)}, \ \
%\la {\bf a}^{(k)},{\bf d}^{(l)} \ra =
%\prod_{i=2}^3 a_{ii}^k d_i^l, \ \
%\la {\bf a},{\bf d} \ra = \prod_{i=2}^3 a_{ii} d_i.
%\la {\bf a}^{(1)},{\bf d}^{(1)} \ra.
%$$
In the non--symmetric case we obtain:
\begin{eqnarray}
2g_0^{(n)}({\bf d}^3)&=&1-s_d-\pi_a+\la {\bf a},{\bf d} \ra\;,
\;\;\;\;\pi_a=\prod_{i=1}^3a_{ii}\;,\label{g0d3n} \\
12g_1^{(n)}({\bf d}^3) & = &
\la {\bf a},{\bf d} \ra
(2\la {\bf a},{\bf d} \ra-3s_d-2\pi_a)
+s_d(s_d+3\pi_a)-\sum_{i\ne j}^3 a_{ii}a_{jj}d_id_j+
\nonumber \\
&&
(d_1d_2+d_1d_3+d_2d_3)+\pi_d-1\;,
\label{g1d3n} \\
12g_2^{(n)}({\bf d}^3)
&=&
\sum_{i=1}^3 A_i \left[(2A_i+1) d_i^3
-\frac{\pi_a}{2}(A_i+2)d_i^2+
\pi_d  d_i + \frac{\pi_d}{2}(2A_i+1)
\right]+
\nonumber \\
&&
\sum_{i,j=1}^3 \left[
C_{ij} d_i^2 d_j
-\frac{\pi_a}{2}B_{ij} d_id_j
-\frac{\pi_d}{2}F_{ij}\right]\;,
\label{g2d3n}
\end{eqnarray}
where %$\la {\bf a}^{(k)},{\bf d}^{(l)} \ra =\prod_{i=1}^3 a_{ii}^k d_i^l$
$$
A_i=a_{ii}-1\;,\;\;
B_{ij}=A_iA_j-A_i-A_j\;,\;\;
C_{ij}=A_i(A_iA_j-A_i-1)\;,\;\;
F_{ij}=A_i(2A_j+1)\;.
$$
We finish the paper giving a compact version of (\ref{gen1a}) through the
Bernoulli polynomials of higher orders
\begin{eqnarray}
g_n({\bf d}^2)=-\frac{B_{n+1}}{n+1} - \frac{1}{d_1d_2(n+1)(n+2)}
[B_{n+2}^{(2)}(0|{\bf d}^2)-B_{n+2}^{(2)}(d_1d_2|{\bf d}^2)]\;,\label{gnd2}
\end{eqnarray}
where
\begin{eqnarray}
B_{n+2}^{(2)}(x|{\bf d}^2)=\sum_{j=0}^{n+2}\sum_{k=0}^j{n+2\choose j}
{j\choose k}d_1^{j-k}d_2^{n+2-j}B_{j-k}B_{n+2-j} x^k\;.
\label{B2expand}
\end{eqnarray}


\begin{thebibliography}{99}
\bibitem{stan96}  R. P. Stanley, {\it Combinatorics and Commutative
                  Algebra}, \\ Birkh\"auser Boston, 2nd ed., (1996)
\bibitem{rods78}  \"O. J. R\"odseth, {\it A note on Brown and Shiue's paper 
                  on a Remark Related to the Frobenius Problem}, 
                  Fibonacci Quaterly, {\bf 32}, 407 (1994)
\bibitem{sylv84}  J. J. Sylvester, {\it Mathematical Questions with Their 
                  Solutions}, \\ Educational Times, {\bf 41}, 171 (1884)
\bibitem{jau93}   T. C. Brown and P. J. Shiue, {\it A Remark Related to the
                  Frobenius Problem}, \\ Fibonacci Quaterly, {\bf 31}, 32 (1993)
\bibitem{fel04}    L. G. Fel, {\it Frobenius Problem for Semigroups ${\sl S}
                   \left(d_1,d_2,d_3\right)$}, \\ preprint, (2004), 
                   [{\tt http://arXiv.org/abs/math.NT/0409331}]
\bibitem{john60}  S. M. Johnson, {\it A Linear Diophantine Problem}, \\
                  Canad. J. Math., {\bf 12}, 390 (1960)
\bibitem{bat53}   H. Bateman and A. Erdel\'yi, {\it Higher Transcendental 
                  Functions}, Vol.1, Chap.1, \\ McGraw-Hill Book Co., NY 
(1953)
\end{thebibliography}
\end{document}